\input AHTOHFIE.STY
\def\Epsilon{{\rm E}}
%\tolerance10000
\hfuzz9pt
%\emergencystretch4mm
%\baselineskip 11pt

%%%%%%%%%%%%%%%%%%%%%%%%%%%%%%%%%%%%%%%%%%%%%%%%%%%%%%%%%%%
%%%%%%%%%%%%%%%%%%%%%%%%%%%%%%%%%%%%%%%%%%%%%%%%%%%%%%%%%%%
%%%%%%%%%%%%%%%%%%%%%%%%%%%%%%%%%%%%%%%%%%%%%%%%%%%%%%%%%%%

\UDC{%
512.543.16 %%Определяющие соотношения
%\hss\hfil\hfil {\sc Preliminary version}
}
\MSC{%
20F06,   %Cancellation theory; application of van Kampen diagrams
20F12,   %Commutator calculus
20F67    %Hyperbolic groups and nonpositively curved groups
}

\title{Commutators
cannot be proper powers in metric small-cancellation torsion-free groups}
\author{%
Elizaveta V. Frenkel$^\flat$
\qquad
%\qquad
%and
%\qquad
\qquad
\qquad
Anton A. Klyachko$^\sharp$
}

\address{
$^\flat$ Faculty
of Further Education,
\quad
$^\sharp$ Faculty of Mechanics and Mathematics 
\\
Moscow State University, Moscow 119991, Leninskie gory
\\
lizzy.frenkel@gmail.com
%\qquad
%\qquad
%\qquad
%\qquad
\qquad
klyachko@mech.math.msu.su
%\qquad
%\quad
}

\footnote{}{\hskip-0.5cm The work of the second author was
supported by the Russian Foundation for Basic Research, project
no.11-01-00945.}

\abstract{%
\narrower 
\narrower 
\narrower 
\narrower 
\narrower 
A nontrivial commutator
cannot be a proper power in a torsion-free group satisfying 
$C'(\lambda)$ small cancellation condition with sufficiently
small~$\lambda$. }

%%%%%%%%%%%%%%%%%%%%%%%%%%%%%%%%%%%%%%%%%%%%%%%%%%%%%%%%%
\s 0. 
Introduction

It is well known that nontrivial commutators cannot be proper powers in 
free groups~[Sch59]. In free products, the situation is more complicated 
but also completely studied [CER94]; some partial results are known about 
amalgamated products [FRR11]. We prove the following fact.

\Th. 
If a torsion-free group satisfies the $C'(\lambda)$ small cancellation 
condition for sufficiently small $\lambda$, then no nontrivial commutator 
can be a proper power in this group.

Recall that a presentation $\pres<X|R>$ satisfies the $C'(\lambda)$ 
\emph{small cancellation condition} if it is symmetrised (i.e. $R$ is 
closed under taking cyclic permutations and inverses) and, for any two 
different relators $r_1, r_2\in R$, the length of their common initial 
segment is less than $\lambda|r_1|$ (see, e.g., [LS80]).

To prove the theorem we use van~Kampen diagrams and the car-crash
lemma [Kl93] (see also [FeR96]); this lemma was applied earlier to
study quite different problems such as equations over groups and 
relative presentations (see, e.g., 
[CG95], 
[FeR96], 
[Kl97],
[FeR98], 
[CG00], 
[CR01], 
[FoR05], 
[Kl05], 
[Kl06a], 
[Kl06b],
[Kl07], 
[Kl09], 
[Le09] 
and 
[KlL12]).

Let us explain our approach on the following toy example. Suppose we want 
to show that a nontrivial commutator cannot be a cube in the free group 
$F(a,b)$. Suppose the contrary. By the Wicks theorem [Wic62], a cyclically 
reduced word is a commutator if and only if 
a cyclic permutation of this word
is graphically equal to 
$xyzx^{-1}y^{-1}z^{-1}$, where $x$, $y$, and $z$ are some reduced words.  
Therefore, we have a graphical equality of the form 
$xyzx^{-1}y^{-1}z^{-1}=www$. This can be described geometrically as 
follows. There is a graph on a torus, and all vertices have degree two, 
except two vertices of degree three, or except one vertex of degree four 
(the latter means that one of the words $x$, $y$, or $z$ is empty). The 
edges of the graph are directed and labeled by letters $a$ and $b$. The 
complement to this graph is homeomorphic to a disk:

\medskip

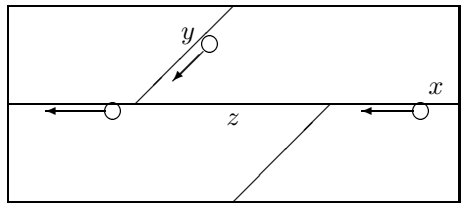
\nobreak
%\vskip1cm plus 1cm minus5mm
%\bigskip
\centerline{\input A.PIC}
%\nobreak%
%\vskip5mm%
%\centerline{Fig. #1}%
%\vskip1cm plus 1cm minus 5mm%
\goodbreak
%\bigskip
%\par
%\noindent

\noindent 
(on this figure, the torus is presented as a rectangle with identified 
opposite sides). Three cars move counterclockwise along the boundary of 
this disk with a constant speed of one edge per minute. This motion is 
periodic with period~$|w|$, i.e. each $|w|$ minutes the cars cyclically 
interchange (each car ``reads'' the word $w$ in $|w|$ minutes). According 
to the car-crash lemma (see~Section~2), every such periodic motion on a 
torus leads to a collision. This collision can occur only at a vertex.  
Indeed, if a car is driving along an edge labeled, say $a$ in the positive 
direction (with respect to the direction of this edge), then all remaining 
cars are also driving edges labeled $a$ in the positive direction at this 
moment (so, they cannot collide). Also, it is easy to see that a collision 
at a vertex contradicts the irreducibility of the boundary label of the 
disk (i.e. of the word $xyzx^{-1}y^{-1}z^{-1}$).

The case of a non-free group is more complicated
due to different structure of the picture: instead of the
simple graph shown above, we have a van~Kampen diagram. However,
for small cancellation groups, this diagram turns out to be very
thin and similar to the graph shown above (see Fig. 2). This allows us
to apply ``automobile technique", although the arguments
become more complicated.

We also utilise the \emph{lowest parameter principle} (as in 
[Ols89]), i.e. we assume that there are fixed small positive numbers  
$ 
\lambda\ll\lambda_1\ll\lambda_2\ll\lambda_3\ll\lambda_4\ll\lambda_5
\ll\lambda_6\ll\lambda_7\ll 1
$ 
and each time an inequality of the form, e.g., 
$\the\year\lambda_2<\lambda_3$ arises, we conclude that it is 
automatically fulfilled by the choice of $\lambda_i$. For reader's 
convenience, we give here a brief list of these parameters.
$$
\settabs
\+\it Parameter\qquad  &Essentially special town\quad  &10\cr
\vbox{
\+\hss\it Parameter\hss&\it Notion\hss                &\S\cr
\+\cr
\+\hss$\lambda$\hss   & $C'(\lambda)$\hss             &0, 1, 3\cr
\+\hss$\lambda_1$\hss & Model\hss                     &1\cr
\+\hss$\lambda_2$\hss & Distance along streets\hss    &1\cr
\+\hss$\lambda_3$\hss & Very near\hss                &4\cr
\+\hss$\lambda_4$\hss & Near\hss                     &4\cr
\+\hss$\lambda_5$\hss & Substantially special town\hss  &3\cr
\+\hss$\lambda_6$\hss & Globally very near\hss       &4\cr
\+\hss$\lambda_7$\hss & Globally near\hss            &4\cr
}
$$
%09.13
The authors thank G. O. Astafurov for useful remarks.
%09.13

%%%%%%%%%%%%%%%%%%%%%%%%%%%%%%%%%%%%%%%%%%%%%%%%%%%%%%%%%%%%
\s 1. 
Small-cancellation maps on a torus with one hole

A \emph{map} on a closed surface is a finite graph on this surface that  
divides the surface into simply connected domains, called \emph{cells} or 
\emph{faces}. Some cells (possibly none) are distinguished and called 
\emph{exterior} or \emph{holes}; the remaining cells are called 
\emph{interior}. The \emph{boundary~$\d\Gamma$} and the \emph{perimeter 
$|\d\Gamma|$} of a cell $\Gamma$ are defined naturally.

A \emph{piece} is a simple path of positive length in the map (i.e. in the 
graph) that joins two vertices (possibly coinciding) of degree other than 
two, and does not pass through other vertices of degree other than two.

We say that a map without vertices of degree one satisfies 
\emph{$C'(\lambda)$ condition} or is a \emph{$C'(\lambda)$-map}, where 
$\lambda$ is a nonnegative real number if the length of any common piece 
of the boundaries of two interior cells $\Gamma_1$ and $\Gamma_2$ is less 
than $\lambda|\d\Gamma_1|$, and the length of any piece that separates an 
interior cell $\Gamma$ and a hole is less than 
$({1\over2}+\lambda)|\d\Gamma|$.

In what follows, we consider only maps on a torus with one hole.

\Lemma 1. 
If $\lambda\ll1$, then the boundary of each interior cell of a 
$C'(\lambda)$-map on the torus with a hole contains at least two pieces 
lying on the hole boundary.

\Proof 
We need also the following simple but useful fact, sometimes called the 
combinatorial Gauss--Bonnet formula.

\proclaim{Weight test \rm[Ger87], [Pri88], see also [MCW02]}.  
If each corner $c$ of a map on a closed surface $S$ is assigned a 
number~$\nu(c)$ {\rm(called the \emph{weight {\rm or the} value of the 
corner $c$})}, then 
$$ 
\sum_v K(v)+\sum_D K(D)+\sum_e K(e)=2\chi(S).  
$$ 
Here, the summations are over all vertices $v$ and all cells $D$ of the 
map, and the values $K(v)$, $K(D)$, and $K(e)$, called the 
\emph{curvatures} of the corresponding vertex, cell, and edge, are defined 
by the formulae
$$
K(v)\:=2-\sum_c \nu(c),
\qquad
K(D)\:=2-\sum_c (1-\nu(c)),
\qquad
K(e)\:=0,
$$
where the first sum is over all corners at the vertex $v$, and the second 
sum is over all corners of the cell $D$.

\noindent 
Let us continue the proof of Lemma 1. To each corner $c$, we assign the 
value $\nu(c)$ by the following rule:
\- 
all corners of the hole have value 1;
\- 
each corner of an interior cell adjacent to a vertex of degree $k+l$ with 
$k$ corners of the hole and $l>0$ corners of interior cells has value 
$\nu(c)=\frac{2-k}{l}$.

\enditem
Thus, the curvature of the hole is 2, and each vertex has nonpositive 
curvature. This curvature is zero if the vertex is adjacent to at least 
one corner of an interior cell. Since the Euler characteristic of the 
torus is 0, the weight test implies that

{\sl the total curvature of all interior cells is at most minus two.}

\noindent
{\bf Step 1. 
The curvature of an interior cell cannot be positive.}

Let $\Gamma$ be a hypothetical interior cell of positive curvature. If the 
boundary of this cell consists of $d$ pieces, then its curvature is at 
most $2-{1\over3}d$ (see Fig. 0).  Hence, the curvature of $\Gamma$ can be 
positive only if $d\le5$.  Such a cell must be adjacent to the hole, 
because otherwise we obtain a contradiction with the small cancellation 
condition. If only one piece of $\Gamma$ lies on the hole boundary, then
$$
0< K(\Gamma)\le 2-{1\over2}-{1\over2}-{1\over3}(d-2)=1-{1\over3}(d-2),
\qqbox{i.e.} d\le 4.
\eqno{(*)}
$$
By definition of a $C'(\lambda)$-map, the length of a common piece of the 
boundaries of $\Gamma$ and the hole does not exceed approximately a half 
of the perimeter of this cell, and the length of each other three or less 
pieces of the boundary is less than $\lambda|\d\Gamma|$, which is 
impossible since $\lambda$ is small enough. If the boundary of $\Gamma$ 
has at least two common pieces with the hole boundary, then the endpoints 
of these pieces of $\Gamma$ have either four corners of weights 
$\le{1\over2}$, two corners of weights $\le{1\over2}$ and one corner of 
weight $\le0$, or two corners of weight $\le0$ (Fig. 0). Thus, 
$K(\Gamma)\le0$.

\goodbreak
%\vskip1cm plus 1cm minus5mm
\bigskip
\centerline{\input 0.PIC}
\nobreak%
%\vskip5mm%
\centerline{Fig. \lowercase{0}}%
%\vskip1cm plus 1cm minus 5mm%
\goodbreak
\bigskip
%\par
%\noindent

\noindent{\bf Step 2. Completion of the proof. At least two pieces of each 
interior cell are adjacent to the hole.}

Indeed, otherwise the cell would have a very big negative curvature (as 
the second inequality of $(*)$ shows, because $d\gg1$ by the small 
cancellation condition), which is a contradiction, because the interior 
cells have nonpositive curvature and the sum of their curvatures is at 
least minus two. This completes the proof of Lemma 1.

\medskip

We call an interior cell \emph{ordinary} if its boundary has two common 
pieces with the hole boundary. An interior cell is called \emph{1-special} 
if its boundary has three common pieces with the hole. An interior cell is 
\emph{2-special} if its boundary has four common pieces with the boundary 
of the hole.

\Lemma 2. 
If $\lambda\ll1$, then each interior cell of a $C'(\lambda)$-map on a 
torus with one hole is either ordinary, 1-special, or 2-special. The 
number of all special cells is at most two; if the map has a 2-special 
cell, then there are no other special cells. All vertices, 
%09.13
except may be two, belong to one of the following classes:
vertices of degree two;
vertices of degree three lying on the hole boundary;
vertices of degree four with exactly two corners of the hole, 
and these corners are not adjacent.
%09.13

\Proof 
Let us use the weight test again, but now we assign the weights to the 
corners in a slightly different manner:

\- a corner of an interior cell at a vertex of degree $>2$
   adjacent to one corner of the hole has value $1\over2$;
\- a corner of an interior cell at a vertex of degree $>2$
   adjacent to two corners of the hole has value 0;
\- all remaining corners have value 1.

\enditem
The curvature of each vertex is nonpositive (if this vertex lies on the 
hole boundary, then, at this vertex, there are one corner of value 1 and 
either two corners of value $1\over2$, or another corner of value 1; if 
this vertex does not belong to the hole boundary, then all corners at this 
vertex have value one). Since all corners of the hole have value one, the 
curvature of the hole is two. The curvature of an interior face is 
nonpositive by Lemma~1. The curvatures of ordinary, 1-special, and 
2-special cell are 0, $-1$, and $-2$, respectively. The curvatures of 
other hypothetical cells (``more special") do not exceed $-3$.

%09.13
The curvature of each vertex that belongs to 
none of the classes listed 
is at most $-1$.
%09.13
But 
according to the weight test, the total curvature must be zero. This 
contradiction completes the proof.

\Lemma 3. 
The edges not lying on the hole boundary form a forest. This forest has at 
most two vertices of degree higher than two.

\Proof 
If some edges lying outside of the hole boundary form a cycle, then, 
cutting the torus along this cycle, one can obtain either 
a sphere with three holes,  
%25-03-13
%либо тор с дырой и сферу с двумя дырами,
a torus with a hole and a sphere with two holes, or  
%либо тор с двумя дырами и сферу с дырой. 
a torus with two holes and a sphere with a hole. 
%25-03-13
Each cell of 
the spherical part shares only a small segment of its boundary with each 
hole, except one, by the small cancellation condition (because the cut was 
done along a cycle consisting of at most four pieces by Lemma 2). 

Contracting the new holes, we obtain a $C'(\lambda_1)$-map on a sphere 
with at most one hole, which is impossible for small~$\lambda_1$. Thus, 
there is no cycles, and we have a forest. The number of vertices of degree 
higher than two in this forest is at most two by Lemma 2, which completes 
the proof.

\medskip

If we contract each edge lying outside of the hole boundary, the surface 
remains a torus by Lemma 3. The obtained map is called the \emph{model} of 
the initial map. The interior cells of the model are called \emph{towns}, 
the vertices having at least three corners of the hole are called 
\emph{junctions}. Pieces of the hole boundary not lying on towns 
boundaries are called \emph{highways}. Each highway connects either two 
towns, two junctions, or a town and a junction. Adding to this model 
``zero length highways'', we obtain a map all whose vertices have degree 
two or three; e.g., if, at some vertex of degree four, there are two 
corners of towns, then we assume that these two towns are connected by a 
highway of zero length. Similarly, we assume that a town and a junction 
are connected by a highway of zero length if this junction lies on the 
boundary of the town; a junction of degree four is a pair of triple 
junctions connected by a highway of zero length.

Thus, a model is a map on the torus, and the cells of this map are the 
hole and towns. By Lemma 2, from each town, two (\emph{ordinary town}), 
three (\emph{1-special town}), or four (\emph{2-special town}) highways 
go; each of these highways leads to another (or the same) town or to a 
junction (of degree three).

Note that, by Lemma 3, the perimeter of each town approximately equals (up 
to $\lambda_1$ multiplied by the perimeter) the perimeter of the 
initial cell, from which this town was obtained by contraction of edges.

\Lemma 4. { 
The model has either  
\item{}
two 1-special towns and no junctions, 
\item{}
one 2-special town and no junctions, 
\item{}
one 1-special town and one triple \(i.e. of degree three\) junction, 
\item{}
or no special towns and two triple junctions \rm (see Fig. 1).  
\enditem 
}

\noindent 
On Figure 1, the torus is represented as a square with identified opposite 
sides. Special towns are labeled by letter~``o''. The dashed lines will be 
explain later.

\Proof 
Let us assign the weights as in the proof of Lemma 2. Then, junctions and 
1-special towns have curvature $-1$; 2-special towns have curvature $-2$; 
the hole has curvature 2; and all remaining towns and vertices have zero 
curvature. So, the assertion of the lemma follows immediately from the 
weight test.

%\vfil\break

%\vskip-30mm

\goodbreak
%\vskip1cm plus 1cm minus5mm
\bigskip
\centerline{\input 1.PIC}
\nobreak%
%\vskip5mm%
\centerline{Fig. \lowercase{1}}%
%\vskip1cm plus 1cm minus 5mm%
\goodbreak
\bigskip
%\par
%\noindent

%\vfil\break

The structure of an initial $C'(\lambda)$-map is shown on Figure 2, where 
the joints of ``chains'' are covered with black circles; there are some 
special cells and junctions behind these circles.

\goodbreak
%\vskip1cm plus 1cm minus5mm
\bigskip
\centerline{\input 2.PIC}
\nobreak%
%\vskip5mm%
\centerline{Fig. \lowercase{2}}%
%\vskip1cm plus 1cm minus 5mm%
\goodbreak
\bigskip
%\par
%\noindent

%%%%%%%%%%%%%%%%%%%%%%%%%%%%%%%%%%%%%%%%%%%%%%%%%%%%%%%%%%%%%%%%%%%%%%%%%%
\s 2. 
Motions

All definitions and facts mentioned in this section are taken from [Kl05].

Consider a map $\Mu$ on a closed oriented surface $S$. A \emph{car} moving 
around a face~$D$ of this map is an orientation preserving covering of the 
boundary $\d D$ of the face $D$ by an oriented circle $R$ (the 
\emph{circle of time}).

Roughly speaking, a car moves along the boundary of its face
counterclockwise (the interior of the face remains on the left
from the car), without U-turns and stops. This motion is periodic.

If the number of cars being at a moment of time $t$ at a point $p$ of the 
1-skeleton of $\Mu$ equals the degree of this point, then we say 
that at the point $p$ at the moment $t$ a \emph{complete collision} 
occurs; the point $p$ is called a \emph{point of complete collision}. 
Points of complete collision lying on edges are called simply \emph{points 
of collision}.

A \emph{multiple motion of period $T$} on a map $\Mu$ is a set of cars 
$\alpha_{D,j}\:R\to\d D$, where $j=1,\dots,d_D$, such that
\item{1)} 
$d_D\ge 1$ (i.e. each face is moved around by at least one car);
\item{2)} 
$\alpha_{D,j}(t+T)=\alpha_{D,j+1}(t)$ for any $t\in R$ and 
$j=\{1,\dots,d_D\}$ (subscripts modulo ${d_D}$, and the addition of points 
of the circle $R$ is defined naturally: $R=\R/l\Z$);
\item{3)} 
there exists a partition of each circle $\partial D$ into $d_D$ arcs (with 
disjoint interiors) such that during the time interval $[0,T]$ each car 
$\alpha_{D,j}$ moves along the $j$-th arc.

\proclaim{Car-crash lemma \rm [Kl05], [Kl97]}. 
For any multiple motion on a map $Mu$ on a closed oriented surface $S$, 
the number of points of complete collision is at least 
$$ 
\chi(S)+\sum_D(d_D-1),
$$
where the sum is over all faces $D$ of $\Mu$.

%%%%%%%%%%%%%%%%%%%%%%%%%%%%%%%%%%%%%%%%%%%%%%%%%%%%%%%%%%%%%%%%%
\s 3. 
Van~Kampen diagrams, buses, and cabs

Let us continue the proof of the main theorem. Suppose that there exists a 
proper power $w^n\ne1$ which is a commutator in a torsion-free 
$C'(\lambda)$-group $G=\pres<X|R>$. The torsion-freeness of a small 
cancellation groups means that all relators in $R$ are not proper powers 
(see, e.g., [LS80]). It is well known (and can be easily proved) that a 
word is a commutator in $G=\pres<X|R>$ if and only if it can be read on 
the hole boundary of some van~Kampen diagram on a torus with one hole. If 
the group satisfies the $C'(\lambda)$ condition, then this van~Kampen 
diagram is a $C'(\lambda)$-map as the following (probably) well-known 
lemma shows. We provide the proof of this lemma for the sake of 
completeness.

\Lemma about powers. 
Suppose that $G=\pres<X|R>$ is a presentation satisfying the $C'(\lambda)$ 
condition, where $\lambda\ll1$, $G$ is a torsion-free group, $w$ is a word 
not conjugate in $G$ to a shorter word, and $n\in\N$. If a word~$v$ is a 
common initial subword of a relator $r\in R$ and $w^n$, then 
$|v|<\({1\over2}+\lambda\)|r|$.

\Proof 
Since $v$ is a subword of $w^n$, it has a form $v=w^kt$, where $t$ is an 
initial segment of $w$. If ${k=0}$, then $|v|\le{1\over2}|r|$, since $w$ 
and its initial subword $t$ are irreducible. Suppose that $k\ge1$. Then 
$w^{k-1}$ is a piece (it is contained twice in $r$). Therefore, the small 
cancellation condition implies that $|w^{k-1}|=(k-1)|w|<\lambda|r|$. If 
$k>1$, we obtain $|w|<\lambda|r|$ and 
$|v|=|w^kt|=|w^{k-1}|+|w|+|t|<3\lambda|r|$.  Hence, the inequality holds 
for all sufficiently small~$\lambda$. If $k=1$, i.e. $v=wt$, then $t$ is a 
piece, and, therefore, $|t|<\lambda|r|$.  Now the required inequality is 
fulfilled, because $|w|\le{1\over2}|r|$ by the irreducibility of $w$.

\medskip

Consider the model (Fig. 1) of a map. On the hole boundary of this model, 
the motion of $n$ cars can be specified in a natural way; these cars are 
called \emph{buses}. Each bus moves with a speed of one edge per minute, 
and, during an $i$-th minute, it drives along an edge labeled by the 
$i$-th letter of $w$ (where $i$ is modulo $|w|$). This motion is periodic 
with period~$|w|$.

Let us draw additional edges called \emph{streets} inside each town in 
such a way that the streets of each town form a tree, 
connecting the exits from the town (i.e. the set of vertices of degree one 
of this tree coincides with the set of endpoints of highways lying on 
the boundary of the town). Moreover, we draw the streets in such a way that 
the distance along streets between any two neighbouring exits from a town 
is approximately equal (up to $\lambda_2$ multiplied by the perimeter of 
the town) to the distance between these exits along the boundary of the 
town. Clearly, this can be done, because each part of the boundary of a 
town is not significantly greater than a half of its perimeter (i.e.  
$\le({1\over2}+\lambda_2)|\d\Gamma|$). The streets and the highways form a 
map (with one cell) on the torus (see Fig. 1 (lower left), where streets 
are drawn as dashed lines). The vertices of degree higher than two 
in special towns are also called \emph{junctions} henceforth. A 
junction (of streets) of degree four is considered as a pair of junctions 
of degree three joined by a highway of zero length. A special 
town~$\Gamma$ is called \emph{substantially special} if all street 
junctions in this town are further than $\lambda_5|\d\Gamma|$ from the 
boundary of the town.

Now, we define a motion of $n$ cars called \emph{cabs} on this map. A 
cab drives along a highway together with the corresponding bus. When 
a cab enters a town, it moves along the streets with a constant 
(approximately unit) speed in such a way that it leaves the town 
simultaneously with the corresponding bus (which drives along the boundary 
of the town). On Figure 3, we draw a special town and positions of a bus 
(black circle) and a cab (white circle) at four moments.

\goodbreak
%\vskip1cm plus 1cm minus5mm
\bigskip
\centerline{\input 3.PIC}
\nobreak%
%\vskip5mm%
\centerline{Fig. \lowercase{3}}%
%\vskip1cm plus 1cm minus 5mm%
\goodbreak
\bigskip
%\par
%\noindent

%%%%%%%%%%%%%%%%%%%%%%%%%%%%%%%%%%%%%%%%%%%%%%%%%%%
\s 4. 
Nearness

Let $x$ and $y$ be points of the graph formed by the highways, streets, 
and boundaries of towns. We say that $x$ is \emph{\(locally\) near $y$} if 
either $x=y$, or $y$ is in a town~$\Gamma$%
\fn{%
Points in zero distance from a town are also considered as being in the 
town.}
and the distance between $x$ and $y$ is at most $\lambda_4|\d\Gamma|$. If 
this distance is at most $\lambda_3|\d\Gamma|$, then we say that $x$ is 
\emph{\(locally\) very near $y$}.

We say that $x$ and $y$ are \emph{globally near} each other if 
the distance between them is at most $\lambda_7|w|$, where $|w|$ is the 
period of the motion. If this distance is at most $\lambda_6|w|$, then 
these points are said to be \emph{globally very near} each other. Note 
that, if $x$ is near $y$, then these points are globally very near each 
other.

%%%%%%%%%%%%%%%%%%%%%%%%%%%%%%%%%%%%%%%%%%%%%%%%%%%%%%%%%%%%%%%%%
\s 5. 
Where do collisions happen?

The car-crash lemma implies that a complete collision of cabs
occurs at least at $n-1$ points. Indeed, the map formed by all
streets and highways has one cell; along the boundary of this cell
$n$ cabs move regularly and the Euler characteristic of the
torus is zero. It remains to understand where can collisions occur.

\Lemma about far-from-junction collisions. 
Suppose that some cabs collide at a point $p$
%09.13
which is not a junction. 
%09.13
Then $p$ is in a town and 
either there is a junction 
%09.13
very 
%09.13
near 
%09.13
to 
%09.13
$p$ or $p$ is on a highway of zero length 
connecting two different special towns with equal labels \rm \(i.e. the 
labels of the corresponding cells of the initial diagram are equal if we 
read them starting from the point $p$\).

\Proof 
Consider several cases.

{\noindent\bf I.}
{\sl Collisions cannot occur on a highway of nonzero length outside  
junctions}.

Indeed, suppose a collision of cabs occurred on an edge labeled by a 
letter $x$. Since cabs move along highways in the same way as buses, this 
means that a bus moves along an edge labeled by $x$ at the moment of 
collision and another bus moves along this edge in the opposite direction, 
i.e. moves along an edge labeled by $x^{-1}$.  This contradicts the 
definition of the motion, because, at each moment of time, all buses drive 
along  edges with the same labels.  A collision at some vertex on a 
highway is also impossible, because, by definition of the motion of buses, 
this would mean that the word $w$ has a subword $xx^{-1}$ that contradicts 
the irreducibility of $w$.

{\smallskip\noindent\bf II.} 
{\sl If two cabs collide in a town $\Gamma$, then at least one of them 
leaves the town in at most $\lambda_1|\d \Gamma|$ minutes after the 
collision.}

The $C'(\lambda_1)$ condition implies that the boundary of $\Gamma$ does 
not have two segments with equal labels of length~$\lambda_1|\d\Gamma|$. 
Therefore, two cabs cannot stay in one town longer 
than~$\lambda_1|\d\Gamma|$ minutes simultaneously.

{\smallskip\noindent\bf III.} 
{\sl If two cabs collide in a town $\Gamma$ and this collision occurred 
further than $\lambda_2|\d\Gamma|$ from the junctions, then, during 
$\lambda_1|\d\Gamma|$ minutes after the collision, precisely one of the 
cabs leaves $\Gamma$.}

If both cabs leave the town during $\lambda_1|\d\Gamma|$ minutes, then the 
distance between two entrances to this town is at most 
$\lambda_2|\d\Gamma|$. This is possible only in a special town near a 
junction (because the distance between two entrances in an ordinary town 
is approximately a half of the perimeter, and there is always a street 
junction between two entrances in a special town).

{\smallskip\noindent\bf IV.} 
{\sl The length of a highway $r$ along which one of the colliding cabs 
leaves the town $\Gamma$ during $\lambda_1|\d\Gamma|$ minutes is at most 
$\lambda_2|\d\Gamma|$.}

Suppose that this highway is longer. Then, one can choose a segment $s$ of 
length at least $\lambda_1|\d\Gamma|$ on this highway such that, before 
the collision, one of the buses is driving along~$s$ while the other one 
is driving along $\d\Gamma$; and, shortly after the collision, one of the 
buses is driving along $s$ (in the opposite direction), while the other 
bus is driving along $\d\Gamma$. This means that the label of the boundary 
of the cell $\Gamma$ contains both $f$ and $f^{-1}$ as subwords, where $f$ 
is the label of $s$.  This is impossible by the small cancellation 
condition.

{\smallskip\noindent\bf V.} 
{\sl The highway $r$ leads to a town $\Delta$ with a label%
\footnote{**$^)$}{%
The label of a town is the label of the corresponding cell of the
initial diagram, i.e. one of the defining relators of~$G$.} 
equals the label of the town $\Gamma$.}

Indeed, shortly after the collision one of the buses passes through the 
highway $r$ and then is driving along the boundary of a town $\Delta$ 
during at least $\lambda_1|\d\Gamma|$ minutes (because we suppose that 
there are no junctions near $p$), while the other bus is driving along the 
boundary of $\Gamma$. The small cancellation condition implies that the 
labels of~$\Gamma$ and $\Delta$ are equal.

{\smallskip\noindent\bf VI.} 
{\sl The highway $r$ has zero length.}

As it was shown above, one of the buses goes from a point $A\in\d\Gamma$ 
to a point $B\in\d\Delta$, while the other bus goes from a point 
$A'\in\d\Delta$ to a point $B'\in\d\Gamma$ (Fig. 4). Each of them is 
driving along a path with the same label $u$ during this time. Thus, the 
labels of $\Gamma$ and $\Delta$ are equal if we read them starting from 
the points $A$ and $A'$; and the same is true, if we read these labels 
starting from points $B'$ and $B$. Therefore, the label~$v$ of the part of 
the boundary of~$\Gamma$ from $B'$ to $A$ (counterclockwise) equals the 
label of the part of the boundary of~$\Delta$ from~$B$ to~$A'$. Hence, 
$(uv)^2=1$ in the free group. This means that $uv=1$. Therefore, the 
highway $r$ has zero length by the irreducibility.

%\vfil\break

\goodbreak
%\vskip1cm plus 1cm minus5mm
\bigskip
\centerline{\input 4.PIC}
\nobreak%
%\vskip5mm%
\centerline{Fig. \lowercase{4}}%
%\vskip1cm plus 1cm minus 5mm%
\goodbreak
\bigskip
%\par
%\noindent

{\smallskip\noindent\bf VII.} 
{\sl The towns $\Delta$ and $\Gamma$ are special.}

Suppose the contrary. The length of each piece of the boundary of an 
ordinary town is approximately the half of the perimeter;  
the pieces of the boundary of a special town also have a big length 
($\ge\lambda_4|\d\Gamma|$), because we assume that there are no junctions 
near the point $p$ (otherwise, there is nothing to prove). The labels of 
towns $\Gamma$ and $\Delta$ are equal, if we read them from the point $p$. 
Therefore, we have a segment of the hole boundary with a label of length 
much more than the length of a defining relator (the part $uv$ on Fig. 5), 
which is impossible by the lemma about powers. This completes the proof.

%\vfil\break

\goodbreak
%\vskip1cm plus 1cm minus5mm
\bigskip
\centerline{\input 5.PIC}
\nobreak%
%\vskip5mm%
\centerline{Fig. \lowercase{5}}%
%\vskip1cm plus 1cm minus 5mm%
\goodbreak
\bigskip
%\par
%\noindent

\Lemma about uniqueness. 
At most one collision can happen globally far from junctions.

\Proof 
By the lemma about far-from-junction collisions, it suffices to prove 
that there cannot be two towns $\Gamma$ and $\Delta$ connected by two 
zero-length highways $p$ and $q$ with collisions happening on each of 
these highways (Fig. 6).

\goodbreak
%\vskip1cm plus 1cm minus5mm
\bigskip
\centerline{\input 6.PIC}
\nobreak%
%\vskip5mm%
\centerline{Fig. \lowercase{6}}%
%\vskip1cm plus 1cm minus 5mm%
\goodbreak
\bigskip
%\par
%\noindent

Three exits from the town $\Gamma$ divides its boundary into three
parts $a$, $b$, and $c$ (listed counterclockwise starting with the
piece between $p$ and $q$). Similar fragments of the boundary of
$\Delta$ are denoted $a'$, $b'$, and $c'$. The labels
of $\Gamma$ and $\Delta$ are equal if are read starting from 
$p$ (or $q$) by lemma about far-from-junction collisions. 
Therefore, the labels of segments $a\cup b'$, $a'\cup b$, $c\cup a'$
and $c'\cup a$ of the hole boundary are equal to subwords
of the label of~$\Gamma$. Hence, by lemma about powers,
these fragments cannot be very long:
$$
|a|+|b'|\pe{1\over2}|\d\Gamma|,
\quad
|a'|+|b|\pe{1\over2}|\d\Gamma|,
\quad
|c|+|a'|\pe{1\over2}|\d\Gamma|,
\quad
|c'|+|a|\pe{1\over2}|\d\Gamma|,
$$
where $x\pe y$ means $x<y+k\lambda_1|\d\Gamma|$ for an
absolute constant $k$. Summing up these four inequalities we obtain
$2|a|+2|a'|+|b|+|b'|+|c|+|c'|\pe 2|\d\Gamma|$. Taking into account
the equality $|a|+|b|+|c|=|a'|+|b'|+|c'|=|\d\Gamma|$, we get
$|a|+|a'|\pe0$ that contradicts the fact that collision points $p$ and $q$ 
are far from a junction (lying in $\Gamma$).

%%%%%%%%%%%%%%%%%%%%%%%%%%%%%%%%%%%%%%%%%%%%%%%%%%%%%%%%%%%%%%%%%
\s 6. 
Neighborhoods of junctions and the proof of the theorem for $n\ne3$

As we already know, there is at most one point of collision lying 
(globally) far from junctions, and there are at most two junctions. 
Therefore, if (globally) near each junction, there was at most one 
collision point, then the total number of collision points would not 
exceed three. By the car-crash lemma, this mean that the multiplicity of 
the motion (i.e. $n$) is at most four which proves the theorem for all 
$n\ge5$. However, a nontrivial commutator in a torsion-free small 
cancellation group cannot be a square (see [Sch80], [Gu89]), and, 
therefore, cannot be a fourth power. Thus, the theorem would be proven for 
all $n$, except 3.

The problem is that we do not know how many collision-points are near 
junctions. The following lemma allows us to overcome this difficulty.

\Lemma about junction neighborhoods. 
In $\lambda_7|w|$-neighborhoods of junctions \(i.e.  globally near them\), 
the schedule of the motion can be modified in such a way that, in each 
connected component of the union of these neighborhoods, at most one 
collision occurs.

\Proof 
Consider the graph formed by the streets and highways; we do not need towns 
anymore. Let us consider neighborhoods of junctions of radius 
$\lambda_7|w|$ in this graph and then thicken them slightly to obtain open 
neighborhoods of junctions on the torus.  The union $U$ of these 
neighborhoods is either two disks, one disk, or one annulus (see the upper 
part of Fig. 7). Now, we remove edges of the graph lying inside $U$ and 
add the boundary of $U$ to this graph; in the annulus case, we add also an 
edge cutting the annulus to make it simply connected (see the lower part 
of Fig. 7).

\goodbreak
%\vskip1cm plus 1cm minus5mm
\bigskip
\centerline{\input 7.PIC}
\nobreak%
%\vskip5mm%
\centerline{Fig. \lowercase{7}}%
%\vskip1cm plus 1cm minus 5mm%
\goodbreak
\bigskip
%\par
%\noindent

We obtain a map on the torus with three (in the first case) or
two (in the second and third cases) faces. Let us specify a motion on
this map. The cabs moving around the old face (the hole)
drive almost as was defined before, but the time they spent inside $U$
they now spend moving along the boundary of $U$. Around
the new faces, new cabs moves as follows. When there are no
other cabs on the boundary of the corresponding face (such a
moment exists, because old cabs spend little time on the boundary
of the domain $U$, less than $\lambda_7$ multiplied by the
period), a new cab drives quickly around almost the entire face, and
then moves slowly along the remaining part of the boundary 
of the face until the end of period. This motion is periodic with the 
same period, the multiplicity of the motion on the new faces is 1 and each 
new cab collides at most once (during the slow motion). This 
completes the proof.

\medskip

To complete the proof of the theorem (for $n\ne3$), it remains to
apply the previous lemma. As it was shown above, there is at most
one point of collision globally far from junctions. Hence, the
argument from the beginning of the section completes the proof
of the theorem for all~$n\ne3$.

%%%%%%%%%%%%%%%%%%%%%%%%%%%%%%%%%%%%%%%%%%%%%%%%%%%%%%%%%%%%%%%%%
\s 7. 
The proof of the theorem for $n=3$

\Lemma about close junctions. 
If the junctions are globally very near each other, then all collisions 
occur globally near the junctions.

\Proof 
Suppose the contrary. By the lemma about far-from-junction collisions, these 
junctions are in different special towns $\Gamma$ and $\Delta$ with equal 
labels, and a collision occurs far from junctions at a point $p$ lying 
on the intersection of the boundaries of these towns (Fig. 8). Moreover, 
the junctions are globally very near each
other, i.e. they are globally very near 
exits from the towns to a (globally very short) highway $r$. In particular, 
this means that the segments $c$ and $c'$ of town boundaries lying in 
opposite sides to the exits to $r$ are equal to about a half of the 
perimeter:% 
\fn{The distance along streets between two exits is approximately equal to 
the distance between these exits along the boundary of the town. 
Therefore, if one of these streets is short, then the sum of lengths of 
the other two streets is approximately a half of the perimeter of the 
town.}

$$
|c|\ge {1\over2}|\d\Gamma|-10\lambda_6|w|\le|c'|.
$$

\goodbreak
%\vskip1cm plus 1cm minus5mm
\bigskip
\centerline{\input 8.PIC}
\nobreak%
%\vskip5mm%
\centerline{Fig. \lowercase{8}}%
%\vskip1cm plus 1cm minus 5mm%
\goodbreak
\bigskip
%\par
%\noindent

If the collision point $p$ is globally far from junctions,
then the segments $a$ and $a'$ of the boundaries of the towns 
connecting $p$ with the endpoints of $r$ are long:
$$
|a'|\ge\lambda_7|w|\le|a|.
$$
Summing up these inequalities we obtain
$$
|c|+|a'|\ge {1\over2}|\d\Gamma|+(\lambda_7-10\lambda_6)|w|\le|c'|+|a|.
$$
However, $c\cup a'$ is a part of the hole boundary and its length cannot 
exceed ${1\over2}|\d\Gamma|+\lambda_1|\d\Gamma|$. This contradiction 
completes the proof.

\Lemma about safe junctions. 
Suppose that a collision occurs at a point $p$; very near $p$, there 
is exactly one junction, and at the moment of the collision, there are 
less than three cabs very near the point $p$. Then the motion of cabs very 
near the point $p$ can be modified in such a way that no collisions in 
this neighborhood of $p$ will happen.

\vfil\break

\goodbreak
%\vskip1cm plus 1cm minus5mm
\bigskip
\centerline{\input 9.PIC}
\nobreak%
%\vskip5mm%
\centerline{Fig. \lowercase{9}}%
%\vskip1cm plus 1cm minus 5mm%
\goodbreak
\bigskip
%\par
%\noindent

\Proof 
Two cabs collide near the junction, and the remaining cabs are far 
from the junction at this moment (Fig. 9, the upper row). Let us slightly 
change the schedule of the motion. Namely, the first cab approaches the 
junction and decelerates a little bit to give to the second cab the 
possibility to pass this junction, then accelerates again and leaves the 
neighborhood of the junction according to its initial schedule 
(Fig. 9, the lower row).

\Lemma about absence of dangerous junctions. 
Suppose that the junctions are not globally very near each other, a 
collision occurs at a point $p$, and there is a junction $j$ very near 
$p$.  Then there cannot be three cabs at the moment of the collision very 
near $p$.

\Proof 
Suppose the contrary. There are five possible cases (Fig. 10).

\vfil\break

\goodbreak
%\vskip1cm plus 1cm minus5mm
\bigskip
\centerline{\input 10.PIC}
\nobreak%
%\vskip5mm%
\centerline{Fig. \lowercase{10}}%
%\vskip1cm plus 1cm minus 5mm%
\goodbreak
\bigskip
%\par
%\noindent

{\noindent\bf I. \sl 
The point $p$ is outside the towns.}

By definition of the nearness, this means that $p=j$, and three buses 
collide at this point simultaneously. Therefore, all three edges leading 
to this junction have the same labels, and the word written on the hole 
boundary is reducible, which is impossible. In what follows, we assume 
that the point $p$ is in a town $\Gamma$. Note that the distance from $p$ 
to the second junction (different from $j$), even if it lies in $\Gamma$, 
is large (at least $\lambda_5|w|$), because these junctions are not 
globally very near each other.

{\smallskip\noindent\bf II. \sl 
The junction $j$ is inside a substantially special town $\Delta$. }

The point $p$ cannot lie in the same town, because this would mean that 
two different long (of length longer than $\lambda_1|\d\Gamma|$) parts of 
the boundary of the town $\Gamma=\Delta$ have equal labels, which is 
impossible by the assumption (similarly to the argument in {\bf II} in the 
proof of the lemma about far-from-junction collisions). Therefore, 
$p\in\Gamma\ne\Delta\ni j$. In this case, repeating word-by-word the 
argument from {\bf V} of the proof of the lemma about far-from-junction 
collisions, one can verify easily that the labels of $\Gamma$ and $\Delta$ 
are equal. This implies that 
$$ 
\lambda_5|\d\Delta|\le\rho(p,j)\le\lambda_3|\d\Gamma|=\lambda_3|\d\Delta| 
\qbox{(from now on $\rho$ means the distance).} 
$$ 
(The first inequality holds, because the town $\Delta$ is substantially 
special; the second inequality holds because $j$ is very near $p$.) This 
is impossible, since $\lambda_3\ll\lambda_5$.

{\smallskip\noindent\bf III. \sl 
There is only one town \(i.e. $\Gamma$\) near the point $p$. }

Consider the $\lambda_4|\d\Gamma|$-neighborhood of $p$; remove the 
$\lambda_3|\d\Gamma|$-neighborhood of this point, and intersect this 
difference with the 1-skeleton of the model. We obtain a disjoint union of 
four paths: $u_1$, $u_2$, $v$, and $w$, where $u_1$ and $u_2$ lie on the 
boundary of $\Gamma$, and $v$ and $w$ lie on highways leading to the 
junction $j$. The label of paths $u_1$, $v$, and $w$ are approximately 
equal (i.e. these three words contain a common subword of the length at 
least $(\lambda_4-10\lambda_3)|\d\Gamma|$), because, along these paths, 
buses are driving simultaneously while approaching the junction. On the 
other hand, the labels of paths $u_2^{-1}$, $v^{-1}$, and $w^{-1}$ are 
also approximately equal in the same sense, because, along these paths, 
buses are moving simultaneously while driving away from the junction.  
Therefore, the labels of paths $u_1$ and $u_2$ are approximately equal 
(contain a common subword of the length at least 
$(\lambda_4-20\lambda_3)|\d\Gamma|$), which contradicts the small 
cancellation condition.

{\smallskip\noindent\bf IV. \sl 
There are precisely two towns near the point $p$. }

In this case, we obtain similarly two long segments with equal
labels on the boundary of the town $\Gamma$, which is a contradiction.

{\smallskip\noindent\bf V. \sl 
There are precisely three towns near the point $p$. }

Let us denote these towns $\Gamma$, $\Delta$, and $\Epsilon$.  Arguing as 
in the two previous cases, we obtain six long (of the length at least 
$(\lambda_4-20\lambda_3)|\d\Gamma|$) paths:  $u_1$, $u_2$, $v_1$, $v_2$, 
$w_1$, and $w_2$. The first two lie on the boundary of $\Gamma$, the 
second two lie on the boundary of $\Delta$, and third two lie on the 
boundary of $\Epsilon$; the labels of $u_1$, $v_1$, and $w_1$ coincide 
(along these paths buses approach the junction simultaneously) and 
labels of $u_2$, $v_2$, and $w_2$ coincide (along this paths buses drive 
away from the junctions simultaneously). By the small cancellation 
condition, this means that the cells of the initial $C'(\lambda)$-map on the 
torus corresponding to the three towns have equal labels, and these labels 
are equal, if one read them starting from paths $u_1$, $v_1$, and $w_1$ or 
starting from paths $u_2$, $v_2$, and $w_2$.

Thus, the paths of boundaries of the three cells between $[u_2^-,u_1^+]$, 
$[v_2^-,v_1^+]$, and $[w_2^-,w_1^+]$ also have a common label $g$, where 
$\pi^-$ and $\pi^+$ denote the origin and the endpoint of a path $\pi$ (we 
assume, that all paths are oriented counterclockwise with respect to the 
hole). Now, consider the following three paths in the model: from $u_1^+$ 
to~$v_2^-$, from~$v_1^+$ to~$w_2^-$, and from $w_1^+$ to~$u_2^-$. They 
also have a common label $h$, because the buses drive along these paths 
simultaneously. Moreover, these six paths form a closed contour in the 
initial diagram, and this contour does not have any cells inside. 
Therefore, $(gh)^3=1$ in the free group, i.e. $h=g^{-1}$.

However, the label of the hole in the initial diagram is a reduced word, 
and, therefore, the boundaries of the three cells have a common point and 
the labels of these three cells are equal if we read them starting from 
this point (Fig~11). This means that the junction $j$ lies outside of 
towns, i.e. only one of these three cells can be special (i.e. contains 
another junction). Arguing as in {\bf VII} in the proof of lemma about 
far-from-junction collisions, we conclude that the label of the hole 
boundary contains a common subword with one of defining relators, and the 
length of this subword exceeds significantly a half of the length of the 
defining relator. If, for example, the towns $\Gamma$ and $\Delta$ on 
Figure~11 are ordinary, then the label of the part of the boundary of the 
hole from $x$ to $y$ contains almost a whole relator as a subword. This 
contradicts the lemma about powers and completes the proof.

\goodbreak
%\vskip1cm plus 1cm minus5mm
\bigskip
\centerline{\input 11.PIC}
\nobreak%
%\vskip5mm%
\centerline{Fig. \lowercase{11}}%
%\vskip1cm plus 1cm minus 5mm%
\goodbreak
\bigskip
%\par
%\noindent

Let us continue the proof of the theorem. By the car-crash lemma, it 
suffices to show that the number of collision points cannot exceed one.

If the junctions are globally very near each other, then, by the lemma 
about close junctions, all points of collision are globally near these 
junctions and, by the lemma about neighborhoods of junctions, we can 
assume that globally near junctions there is at most one point of 
collision. This means that there is at most one collision point, as 
required.

If the junctions are not globally very near each other, then, applying the 
lemma about absence of dangerous junctions and the lemma about safe 
junctions, we again obtain a motion with at most one collision point. This 
completes the proof of the theorem.

%%%%%%%%%%%%%%%%%%%%%%%%%%%%%%%%%%%%%%%%%%%%%%%%%%%%%%%%%%%%%%%%%
\REFERENCES

%\[Kl05]
%Klyachko Anton. A.
%Conjecture Kervaire--Laudenbach and presentation simple groups
%{// algebra and logic.} 2005. {Th. 44}. {\number4}. S. 399--437.
%See also
%{arXiv:math.GR/0409146}
\[Kl05]
Klyachko Ant. A.
The Kervaire--Laudenbach conjecture and presentations of simple groups
{// Algebra i Logika}. 2005. {T. 44}. {no.4}. P. 399--437.
See also
{arXiv:math.GR/0409146}
%\[Kl05]
%Klyachko Ant. A.
%The Kervaire--Laudenbach conjecture and presentations of simple groups
%{// Algebra and Logic.} 2005. {T. 44}. {no.4}. P. 219--242.
%See also
%{arXiv:math.GR/0409146}.

%\[Kl06a]
%Klyachko Anton. A.
%As generalize known results about equations over groups
%{// Mat. Zametki}. 2006. {Th.79}. \number3. S. 409--419.
%See also
%arXiv:math.GR/0406382.
\[Kl06a]
Klyachko Ant. A.
How to generalize known results on equations over groups
{// Mat. Zametki}. 2006. {T.79}. no.3. P.409--419.
See also
arXiv:math.GR/0406382.

%\[Kl06b]
%Klyachko Anton. A.
%SQ-universality relative presentations with one the relation
%{// Mat. collection}. 2006. {T.197}. \number 10. S.87--108.
%See also
%arXiv:math.GR/0603468.
\[Kl06b]
Klyachko Ant. A.
The SQ-universality of one-relator relative presentations
{// Mat. Sbornik}. 2006. {T.197}. no.10. P.87--108.
See also
arXiv:math.GR/0603468.

%\[Kl07]
%Klyachko Anton. A.
%Free subgroup relative presentations with one the relation
%{// algebra and logic}. 2007. {Th.46}. \number 3. S.290--298.
%See also
%arXiv:math.GR/0510582.
\[Kl06b]
Klyachko Ant. A.
Free subgroups of one-relator relative presentations
{// Algebra i Logika}. 2007. {V.46}. no.3. P.290--298
See also
arXiv:math.GR/0510582.

%\[LS80]
%Lyndon R., Schupp P.
%Combinatorial theory groups.
%Moscow: Mir, 1980.
\[LS77]
Lyndon R.C., Schupp P.E.
{Combinatorial Group Theory},
Springer-Verlag, Berlin/Heidelberg/New~York, 1977.

\[Ols89]
Olshanskii A. Yu. 
Geometry of defining relations in groups.  Moscow: Nauka, 1989.

\[CG95]
Clifford A., Goldstein R.Z.
Tesselations of $S^2$ and equations over torsion-free groups
{// Proc. Edinburgh Math. Soc.} 1995. {V.38}. P.485--493.

\[CG00]
Clifford A., Goldstein R.Z.
Equations with torsion-free coefficients
{// Proc. Edinburgh Math. Soc.} 2000. {V.43}. P.295--307.

\[CR01]
Cohen M. M., Rourke C.
The surjectivity problem for one-generator, one-relator extensions of
torsion-free groups
{// Geometry \& Topology}. 2001. {V.5}. P.127--142.
See also
arXiv:math.GR/0009101

\[CER94]
Comerford L. P.,  Edmunds C. C.,  Rosenberger G.
Commutators as powers in free products of groups.
Proc. Amer. Math. Soc. 122 (1994), 47-52.

\[FeR96]
Fenn R., Rourke C.
Klyachko's methods and the solution of equations over torsion-free groups
{// L'Enseignment Math\'ematique.} 1996. {T.42}. P.49--74.

\[FeR98]
Fenn R., Rourke C.
Characterisation of a class of equations with solution over torsion-free
groups,
from {``The Epstein Birthday Schrift"},
{(I. Rivin, C. Rourke and C. Series, editors)},
{Geometry and Topology Monographs.} 1998. {V.1}. P.159-166.

\[FRR11]
Fine B., Rosenberger A., Rosenberger, G.
Quadratic properties in group amalgams
//
Journal of Group Theory. V.14, no.5, P.657-671

\[FoR05]
Forester M., Rourke C.
Diagrams and the second homotopy group
{// Comm. Anal. Geom.} 2005. {V.13}. P.801-820.
See also
arXiv:math.AT/0306088

\[Ger87]
Gersten S.M.
Reducible diagrams and equations over groups.
In {Essays in group theory}, P.15--73.
Springer, New York-Berlin, 1987.

\[Gu89]
Zhi-Bin Gu.
Hyperbolic surfaces and quadratic equations in groups
//
Proc. Amer. Math. Soc. 107 (1989), 859-866.

\[Kl93]
Klyachko Ant. A.
A funny property of a sphere and equations over groups
{// Comm. Algebra}. 1993. {V.21}. P.2555--2575.

\[Kl97]
Klyachko Ant. A.
Asphericity tests
{// IJAC}. 1997. {V.7}. P.415--431.

\[Kl09]
Klyachko Ant. A.
The structure of one-relator relative presentations and their centres
//
Journal of Group Theory, 2009, 12:6, 923--947.
See also
arXiv:math.GR/0701308

\[KlL12]
Klyachko Ant. A., Lurye D. E.
Relative hyperbolicity and similar
properties of one-generator one-relator relative presentations with
powered unimodular relator
//
J. Pure Appl. Algebra,
216:3, (2012), 524-534.
See also arXiv:1010.4220.

\[Le09]
Le Thi Giang.
The relative hyperbolicity of one-relator relative presentations
//
Journal of Group Theory. 2009. 12:6, 949--959.
See also
arXiv:0807.2487

\[MCW02]
McCammond J.P., Wise D.T.
Fans and ladders in small cancellation theory
//
{Proc. London Math. Soc. (3)}. 2002. 84(3):599--644.

\[Pri88]
Pride S.J.
Star-complexes, and the dependence problems for hyperbolic complexes
//
{Glasgow Math. J.} 1988. 30(2):155--170.

\[Sch59]
Sch\"utzenberger M. P.  Sur  l'equation $a^{2+n}=b^{2+m}c^{2+p}$ dans
un groupe libre
//
C. R. Acad. Sci.  Paris S\'er. I Math. 248 (1959), 2435--2436.

\[Sch80]
Schupp P .E.
Quadratic equations in groups, cancellation diagrams on compact surfaces,
and automorphisms of surface groups.
Word problems II. Amsterdam:North-Holland, 1980.

\[Wic62]
Wicks M. J.
Commutators in free products
//
J. London Math. Soc. 37 (1962), 433--444.

\end